\documentstyle[12pt]{article}

\newcommand{\const}{\mathop{\rm const}\limits}

\newcommand{\mod}{\mathop{\rm mod}\limits}

\newcommand{\vraisup}{\mathop{\rm vraisup}\limits}

\newcommand{\supp}{\mathop{\rm supp}\limits}

\begin{document}
\begin{center}

{\bf BILATERAL SMALL LEBESGUE SPACES.}\\

\vspace{3mm}

{\sc Eugene Ostrovsky, Leonid Sirota}\\

\vspace{3mm}

 {\it Department of Mathematic, HADAS company, \\
56209, Rosh Ha-Ain, Ha Melecha  street, 22, ISRAEL; \\
 e-mail: galo@list.ru }\\
{\it Bar-Ilan University, 59200, Ramat Gan, ISRAEL; \\
 e-mail: sirota@zahav.net.il }\\

\vspace{3mm}

              Abstract.\\
\end{center}

\textwidth = 8 cm
\begin{verse}
\normalsize
 \hspace{5mm} In this article we investigate the so-called
 \hspace{5mm}  Bilateral Small Lebesgue Spaces: prove that they
 \hspace{5mm} are associated to the  Grand Lebesgue spaces,
 \hspace{5mm} calculate its fundamental functions and Boyd's indices,
 \hspace{5mm} find its dual spaces etc.\\

\vspace{3mm}

Key words and phrases: {\it Grand and Small Lebesgue spaces,
rearrangement invariant (r.i.) spaces, fundamental functions,
Boyd's index, associate and dual spaces.}\\

\vspace{3mm}

{\it Mathematics Subject Classification.} Primary (1991) 37B30,
33K55, Secondary (2000) 34A34, 65M20, 42B25. \\

\end{verse}

\vspace{3mm}

{\bf 1. \ Definitions. Preliminaries. Grand Lebesgue Spaces.}\\

\vspace{3mm}

 Let $ (X, \Sigma,\mu) $ be a measurable space with
non - trivial measure $ \mu: \ \exists A \in \Sigma, \mu(A) \in
(0,\mu(X)). $ We will consider in this article the following conditions on
the measure $ \mu $ or more exactly on the measurable space
$ (X, \Sigma,\mu): $ \par

{\bf 1.}  {\it Finiteness, or the probabilistic case:} $  \mu(X) = 1. $ \\
It is easy to understand that the case $ \mu(X) \in (0,\infty) $ can be
considered by means of renorming of measure. \par

{\bf 2.} {\it Sigma-finiteness:}  There exists the sequence $ \{A(n), \
A(n) \in \Sigma,  \} $ such that

$$
\mu(A(n)) \in (0, \infty), \ \cup_{n=1}^{\infty} A(n) = X.
$$

{\bf 3.} {\it Diffuseness:}

$$
 \forall A \in \Sigma, 0 < \mu(A) < \infty \
\exists B \subset A, \mu(B) = \mu(A)/2.
$$

{\bf 4.} {\it Separability: } The metric space $ (\Sigma, \rho), $ where the
distance $ \rho $ is defined as usually as

$$
\rho(A,B) = \mu(A \setminus B) + \mu(B \setminus A), \ A,B \in \Sigma
$$
is separable.\\

{\bf 5.} {\it Resonant property.} This means that the measure $ \mu $
is nonatomic  or the set $ X $ consists only  on the countable (finite case
is trivial) set of points with equal non \ - \ zero finite measure. \par

 More information about these definitions see in the classical monograph
M.M.Rao \cite{Rao3}.\par

\vspace{2mm}

 {\it We suppose in this report the measure $ \mu $ is sigma-finite. } \par

\vspace{2mm}

 Define as usually for arbitrary measurable function $ f: X \to R^1 \ $
$$
|f|_p = \left(\int_X |f(x)|^p \ \mu(dx) \right)^{1/p}, \ p \ge 1;
$$
$ L_p = L(p) = L(p; X,\mu) = \{f, |f|_p < \infty \}. $
 Let $ a = \const \ge 1, b = \const \in (a,\infty], $ and let $ \psi = \psi(p) $
be some positive continuous on the {\it open} interval $ (a,b) $
function, such that there exists a measurable function $ f: X
\to R $ for which
$$
f(\cdot) \in \cap_{p \in (a,b) } L_p, \ \psi(p) = |f|_p, \ p \in
(a,b). \eqno(1.0)
$$
  We will say that the equality (1.0) and the function $ f(\cdot) $ from
(1.0) is the {\it representation } of the function $ \psi. $ \par
 The set of all those functions we will denote $ \Psi: \ \Psi = \Psi(a,b) =
\{ \psi(\cdot) \}. $ The complete description of this functions
see, for example, in  \cite{Ostrovsky2}, p. 21-27,  \cite{Ostrovsky3}. \par

 We extend the set $ \Psi $ as follows: for arbitrary $ \psi(\cdot) \in
\Psi \ $ define
$$
 EX \Psi \stackrel{def}{=} EX \Psi(a,b) = \{ \nu = \nu(p) \} =
$$
$$
\{ \nu: \exists \psi(\cdot) \in \Psi: 0 < \inf_{p \in (a,b)}
\psi(p)/\nu(p) \le \sup_{p \in (a,b)} \psi(p)/\nu(p) < \infty \},
$$
 $$
U\Psi \stackrel{def}{=} U \Psi(a,b) = \{\psi = \psi(p),  \
\inf_{p \in (a,b) } \psi(p) > 0. \}
$$
and the function $ p \to \psi(p), \ p \in (a,b) $ is continuous.
\par
 We define in the case $ b = \infty \ \psi(b-0) = \lim_{p \to \infty}
\psi(p). $ \par

\vspace{2mm}

{\bf Definition 1.} Let $ \psi(\cdot) \in
U\Psi(a,b). $ The space $ BSGL(\psi) = G(\psi) = G(X,\psi) =
G(X,\psi, \mu) = G(X,\psi,\mu, a,b) $ (Bilateral Grand Lebesgue
Space) consist on all the measurable functions $ f: X \to R $ with
finite norm
$$
||f||G(\psi) \stackrel{def}{=} \sup_{p \in (a,b)} \left[
|f|_p/\psi(p) \right].
$$

 We denote $ B (\psi) = \{p: \ |\psi(p)| < \infty \} \ $ the bounded part
of a function $ f.$\par
 Note that if $ \max(\psi(a+0), \psi(b-0)) < \infty, b < \infty, $ then the
space $ G(\psi) $ coincides with the direct sum $ L_a + L_b, $ and
if $ \min(\psi(a+0), \psi(b-0)) = 0, $ then $ G(\psi) = \{0 \}.$
\par

\vspace{2mm}

 {\sc Hence, we can and will assume further that }

$$
 \max(\psi(a+0), \psi(b-0)) = \infty, \ \inf_{p \in [a,b]} \psi(p) > 0.
$$

\vspace{2mm}

 In the considered case the spaces $ G(\psi), \ \psi \in U\Psi $ are
non-trivial: arbitrary bounded: $ \vraisup_x |f(x)| < \infty $
measurable function $ f: X \to R $ with finite support:

$$
\mu(\supp \ (f) ) < \infty), \ \supp(f) \stackrel{def}{=} \{x, |f(x)| > 0 \} \
$$
 belongs to arbitrary space $ G(\psi) \ \forall \psi \in U \Psi. $ \par

\vspace{3mm}

{\sc The Bilateral Small Lebesgue space  $ BSL(\psi) = SL(\psi) = S(\psi) $
may be defined as an associate space to the r.i.  Grand Lebesgue space }

$$
 BSL(\psi) = S(\psi) = SL(\psi) \stackrel{def}{=} G^/(\psi). \eqno(1.1)
$$

\vspace{3mm}

 Recall that  the associate space to the Banach functional space $ B $ over
introduced measurable space is defined as a (closed) linear subspace of dual
(or conjugate) space $ B^*, $  consisting on all the functionals of a view:

$$
l(f) = l_g(f) = \int_X f(x) \ g(x) \ d \mu, \ f \in B, \ g \in B^/. \eqno(1.2)
$$
 We may say that the function $ g $ generates the  functional $ l_g $ and will
identify as ordinary the functional $ l_g $ and the function $ g.$ \par

\vspace{3mm}

{\bf We investigate in this paper some properties of the Bilateral Small Lebesgue $ BSL $ spaces: find its explicit view, give some examples of norm
calculations  and estimations for the functions belonging to this spaces,
prove that they  obeys the Fatou, Lebesgue and Absolute continuous norm
property,
state its  separability and non-reflexivity, find its fundamental functions and Boyd's indices,  conditions for convergence etc.}\par

\vspace{3mm}	

 The Small Lebesgue spaces were intensively studied in the last ten years,
see e.g.,\cite{Capone1}, \cite{Fiorenza1},  \cite{Fiorenza2},
\cite{Fiorenza3}, \cite{Fiorenza4}, \cite{Iwaniec1}, \cite{Iwaniec2}
 etc. \par
 We intend to continue this investigations: consider some examples, calculate
Boyd's indices etc., and generalize to the case of {\it bilateral
spaces}.\par
 This means that we can consider the case when $ p \in (a,b), $ where $ a $
may differ from $ 1 $  and $ b $ may differ from $ \infty, $
i.e. we may consider not only the cases $ p \in [1, p_0) $ or $ p \in (p_0, \infty). $ \par

For example, the space $ L^{b)}, b \in [1, \infty) $
(in the notations of the articles \cite{Capone1}, \cite{Fiorenza1},  \cite{Fiorenza2}, \cite{Fiorenza3}) coincides with our space
$ G(1,b,0,1/b). $ \par

 All the properties of the BSL spaces which have a proof alike  one in the
case one \ - \ side  spaces,  for instance, provided in the articles
\cite{Capone1}, \cite{Fiorenza1},  \cite{Fiorenza2},
\cite{Fiorenza3}, \cite{Fiorenza4}, \cite{Iwaniec1}, \cite{Iwaniec2}, will be
described  here very briefly.\par

 Note that the $ G(\psi) $ spaces are the particular case of interpolation spaces (so-called $ \Sigma - $ spaces) \cite{Carro1}, \cite{Davis1},
 \cite{Jawerth1}, \cite{Karadzhov1},  \cite{Steigenwalt1}. But we hope that the
our direct representation of these spaces (definition 3) \ is more
convenient for investigation and application.\par
 In the {\it probabilistic} case $ \mu(X) = 1 $ the spaces
$\ G(\psi) \ $ spaces appeared
In the
article \cite{Kozatchenko1}, where are applied  to the theory of random fields. \par

Now
 we consider a very important for further considerations the examples
of $ G(\psi) $ spaces. Let $ a = \const \ge 1, b = \const \in (a,
\infty]; \alpha, \beta = \const. $ Assume also that at $ b < \infty
\ \min(\alpha,\beta) \ge 0 $ and denote by $ h $ the (unique) root
of equation
$$
(h-a)^{\alpha} = (b-h)^{\beta}, \ a < h < b;
  \ \zeta(p) = \zeta(a,b; \alpha,\beta; p) =
$$
$$
(p-a)^{\alpha}, \ p \in (a,h); \ \zeta(a,b; \alpha,\beta;p) =
(b-p)^{\beta}, \ p \in [h,b);
$$
and in the case $ b = \infty $ assume that
 $ \alpha \ge 0, \beta < 0; $ denote
by $ h $ the (unique) root of equation
 $ (h-a)^{ \alpha} = h^{ \beta}, h > a; $ define in this case

$$
\zeta(p) = \zeta(a,b;\alpha,\beta;p) = (p-a)^{\alpha}, \ p \in
(a,h); \ p \ge h \ \Rightarrow \zeta(p) = p^{\beta}.
$$
Note that at
 $ b = \infty \ \Rightarrow \zeta(p) \asymp (p-a)^{\alpha} \
p^{-\alpha + \beta} \asymp \min \{(p-a)^{\alpha}, p^{\beta} \}, \
p \in (a,\infty) $ and that at
 $ b < \infty \ \Rightarrow \zeta(p) \asymp
(p-a)^{\alpha} (b-p)^{\beta} \asymp \min \{(p-a)^{\alpha},
(b-p)^{\beta} \}, \ p \in (a,b). $
 In the case $ \alpha = 0, b < \infty $ we define $ \zeta(p) = (b-p)^{\beta},
 \ p \in (a,b); $ analogously, if $ \beta = 0, \ b < \infty \ \zeta(p) =
(p-a)^{\alpha}, \ p \in (a,b). $ \par
  We will denote also by the symbols $ C_j, j \ge 1 $ some "constructive" finite
non-essentially positive constants which does not depend on the
$ p,n,x $ etc. By definition, the indicator function of a measurable set
$ A, \ A \in \Sigma $ may be defined as usually:

$$
I(A) = I(A,x) = I(x \in A) = 1, x \in A; \ I(A) = 0, x \notin A.
$$

\vspace{3mm}

{\bf Definition 2.} \\

\vspace{3mm}

The space $ G = G_X = G_X(a,b;\alpha,\beta)=
G(a,b) = G(a,b; \alpha,\beta) $ consists on all measurable
functions
 $ f: X \to R^1 $ with finite norm
$$
||f||G(a,b; \alpha,\beta) = \sup_{p \in (a,b)} \left[ |f|_p \cdot
\zeta(a,b; \alpha,\beta;p) \right]. \eqno(1.3)
$$
 On the other word, the space $ G_X(a,b;\alpha,\beta) $ is the particular case of $ G(\psi) $ space4s with the $ \psi $ function of a view: $ \psi(p) =
1/\zeta(p),  $ and $ p \in (a,b). $ \par

 Notice that if $ \psi \in \Psi, \ p \in (a,b), \ b < \infty, $ and
$$
\psi(p) \sim (p \ - \ a)^{-\alpha}, p \to a + 0; \ \psi(p) \sim
(b \ - \ p)^{-\beta},
$$
$ p \to b \ – \ 0; \ \alpha, \beta \ge 0, $ then the space $ G(\psi, a,b) $ is equivalent to the space $ G(a,b; \alpha,\beta). $ \par

\vspace{3mm}

 {\bf Corollary 2. } As long as the cases $ \alpha \le 0; \ b <
\infty, \beta \le 0 $ and $ \ b = \infty, \beta \ge 0 $ are
trivial, we will assume further that either $ 1 \le a < b <
\infty, \min(\alpha,\beta) > 0, $ or $ 1 \le a, b = \infty, \alpha
\ge 0, \beta < 0. $ \par

 The introduced $ G(\psi) $  spaces are some generalization of the
so-called {\it Grand Lebesgue spaces}, see \cite{Ostrovsky1}, \cite{Ostrovsky3}.\par

 The complete description of the spaces conjugated to the (linear topological)
spaces $ \cap_p L_p, $ see in \cite{Davis1},  \cite{Steigenwalt1}.
 The spaces which are conjugate
to Orlicz's spaces are described, e.g., in \cite{Rao1}, pp. 128-135, \cite{Rao2}, chapter 3. \par

 We denote by $ G^o = G^o_X(\psi), \ \psi \in U\Psi $ the closed subspace of
 $ G(\psi), $ consisting on all the functions $ f, $ satisfying the following condition:
$$
\lim_{\psi(p) \to \infty} |f|_p/\psi(p) = 0;
$$
 and denote by $ GB = GB(\psi) $ the closed span in the norm $ G(\psi) $
the set of all bounded: $ \vraisup_x |f(x)| < \infty $
 measurable functions with finite support: $ \mu(\supp \ |f|) < \infty. $ \par

\vspace{2mm}

 Another definition: for a two functions
$ \nu_1(\cdot), \ \nu_2(\cdot) \in U\Psi $ we will write $ \nu_1
<< \nu_2, $ iff
$$
\lim_{\nu_2(p) \to \infty} \nu_1(p)/\nu_2(p) = 0.
$$
 If for some $ \nu_1(\cdot), \nu_2(\cdot) \in U\Psi, \ \nu_1 <<
\nu_2 $ and $||f||G(\nu_1) < \infty, $ then $ f \in G^0(\nu_2).$
\par

\vspace{2mm}

{\sc Examples.} \\

\vspace{2mm}

  We consider now some important examples, which are some generalizations
of considered one in the works \cite{Ostrovsky1},  \cite{Ostrovsky3}.\par

 Let $ X = R, \ \mu(dx) = dx, 1 \le a < b < \infty, \gamma = \const
> -1/a, \ \nu = \const >-1/b, \ p \in (a,b), $
$$
f_{a,\gamma} = f_{a,\gamma}(x) = I(|x| \ge 1) \cdot |x|^{-1/a}
(|\log |x| \ |)^{\gamma},
$$	
$$
g_{b,\nu} = g_{b,\nu}(x) = I(|x| < 1) \cdot |x|^{-1/b} |\log
|x||^{\nu},
$$
$$
h_m(x) = (\log |x|)^{1/m} I(|x|<1), \ m = \const > 0,
$$
$$
 f_{a,b;\gamma,\nu}(x) = f_{a,\gamma}(x) + g_{b,\nu}(x), \
 g_{a,\gamma,m}(x) = h_m(x) + f_{a,\gamma}(x),
$$
$$
\psi^p_{a,b;\gamma,\nu}(p) = 2(1-p/b)^{-p\nu-1} \ \Gamma(p \gamma
+ 1) + 2 (p/a - 1)^{-p \gamma-1} \Gamma(p \nu + 1),
$$
$$
\psi^p_{a,\gamma,m}(x) = 2(p/a-1)^{-p\gamma-1} \Gamma(p \gamma+1)
+ 2 \Gamma((p/m) + 1),
$$
$ \Gamma(\cdot) $ is usually Gamma-function. \par
 We find by the direct calculation:
$$
\left| f_{a,b;\gamma,\nu} \right|^p_p =
\psi^p_{a,b;\gamma,\nu}(p); \ \left|g_{a,\gamma,m} \right|^p_p =
\psi^p_{a,\gamma,m}(p).
$$
 Therefore,
$$
\psi_{a,b;\gamma,\nu}(\cdot) \in \Psi(a,b), \
\psi_{a,\gamma,m}(\cdot) \in \Psi(a,\infty).
$$
Further,
$$
f_{a,b;\gamma,\nu}(\cdot) \in G(a,b; \gamma+ 1/a, \nu +1/b)
\setminus G^o(a,b; \gamma+ 1/a, \nu + 1/b),
$$
$$
g_{a,\gamma,m}(\cdot) \in G \setminus G^0(a,\infty; \gamma+ 1/a,
-1/m),
$$
and $ \forall \Delta \in (0, 1) \ \Rightarrow f_{a,b,\alpha,\beta}
\notin $
$$
 G(a,b;(1-\Delta)(\gamma+1/a),
\nu+1/b)) \cup G(a,b;1/a, (1-\Delta)(\nu+1/b),
$$
$$
g_{a,\gamma,m}(\cdot) \in G \setminus G^o(a,\infty; \gamma + 1/a;
- 1/m).
$$
 More generally, let us consider the following examples. Let $ L = L(z),
\ z \in (0, \infty) $ be slowly varying as $ z \to \infty $ continuous
positive function. \par

The reader can receive  more information  about the slowly varying function
in the monograph \cite{Seneta1}. \par

 Denote for the function $ \psi =
\psi(a,b; \alpha, \beta; p); \ p \in (a,b) $

$$
\psi_{La}(a,b; \alpha, \beta; p) = \psi(a,b; \alpha,\beta) \ L(a/(p-a));
$$

$$
\psi_{Lb}(a,b; \alpha, \beta; p) = \psi(a,b; \alpha,\beta) \ L(b/(b-p));
$$

$$
\psi_{La,Lb}(a,b; \alpha, \beta; p) = \max(\psi_{La}(a,b;\alpha, \beta; p),
\psi_{Lb}(a,b;\alpha,\beta; p) ).
$$
 We consider here the case $ X = R^n $ with usually norm for the $ n \ - $
dimensional vector
$ \vec{x} =  x = (x_1, x_2, \ldots, x_n) \in X: \ |x| = (\sum_{i=1}^n x^2_i)^{1/2} $ equipped
with the (weight) measure
$$
\mu_{\sigma}(A) = \int_A |x|^{\sigma} \ dx_1 dx_2 \ldots dx_n, \
\sigma = \const.     \eqno(1.4)
$$
 Define the function

$$
f_L(x) = I(|x| < 1) \ |x|^{-1/b} |\log |x||^{\gamma} \ L(|\log |x||);
$$

$ b = \const, \ (n+\sigma)b > 1, \ p \in [1, (n+\sigma)b). $ We get after
some calculations using multidimensional polar coordinates and well-known
properties of slowly varying functions [28, p. 30-44 ]:

$$
f_L(\cdot) \in G \setminus G^o(\psi_{Lb})(a,b(n+\sigma); 0, \gamma + 1/b; p)
 \ L(b/(b(n + \sigma) – p)),
$$
$ 1 \le a < b(n+\sigma). $ \par

We define analogously the function 	

$$
g_L(x) = I(|x| > 1) \ |x|^{-1/a} |\log |x||^{\gamma} \ L(|\log |x||);
$$

$ a = \const, \ (n+\sigma)a \ge 1, \ p \in (a(n+\sigma), b), \ b \in
(a(n + \sigma), \infty). $ We receive:

$$
g_L(\cdot) \in G \setminus G^o(\psi_{La})(a(n+\sigma),b; \gamma + 1/a,0; p)
 \ L(a/(p - a(n + \sigma))).
$$
 Correspondingly, if $ 1 \le a(n+\sigma) < b(n+\sigma) < \infty, $ then

$$
f_L + g_L \in G \setminus G^o(\psi_{La,Lb}).
$$
 Let now
$$
\omega(n) = \pi^{n/2}/\Gamma(n/2 + 1), \ \Omega(n) = n \omega(n) = 2 \pi^{n/2}
/\Gamma(n/2),
$$

$$
R = R(\sigma,n) = [(\sigma + n)/\Omega(n)]^{1/(\sigma + n)}, \ \sigma + n > 0,
$$
such that
$$
\mu_{\sigma} \{x: \ |x| < R \} = 1,
$$
and let $ h = h(|x|) $ be some non–negative measurable function, $ h(x) = 0 $ if
$ |x| \ge R(\sigma,n); \ u \ge \exp(2) \ \Rightarrow $

$$
\mu_{\sigma} \{x: h(|x|) > u) \} = \min(1, \exp \left(- W(\log u)) \right),
$$
where $ W = W(z) $ is twice differentiable strong convex in the domain
$ z \in [2, \infty) $ strong increasing function. Denote by

$$
W^*(p) = \sup_{z > 2}(pz-W(z))
$$
the Young – Fenchel transform of the function $ W(\cdot), $ and define the function
$$
\psi(p) = \exp \left( W^*(p)/p \right).
$$
 It follows from the theory of Orlicz's spaces
[22, p. 12 - 18] that at $ p \in [1, \infty) $

$$
|h|_p \asymp \psi(p).
$$

 Another examples. Put for $ X = R^1, \ \sigma = 0, $
$$
f^{(a,b; \alpha,\beta)}(x) = |x|^{-1/b} \exp \left(C_1|\log |x| \
|^{1-\alpha} \right) I(|x| < 1) +
$$
$$
I(|x|\ge 1) \ |x|^{1/a} \ \exp \left(C_2 (\log |x|)^{1 - \beta}
\right);
$$
$ 1 \le a < b < \infty; \alpha,\beta = \const \in (0,1). $ We
obtain by direct computation using the saddle \ - \ point method:

$$
\log \left|f^{(a,b;\alpha,\beta)}(\cdot) \right|_p \asymp
(p-a)^{1-1/\alpha} + (b-p)^{1-1/\beta}, \ p \in (a,b).
$$

 {\it It is known that the spaces } $ G(\psi) $ {\it with respect
to the ordinary operations and introduced norm } $
||\cdot||G(\psi) $ {\it are Banach functional Spaces  in the terminology of
a book  \cite{Bennet1}, they obey the Fatou property etc.} \par

\vspace{2mm}
{\sc Properties. Fundamental function.}\\
\vspace{2mm}

 See also \cite{Fiorenza1}, \cite{Fiorenza2}, \cite{Capone1},
\cite{Ostrovsky1}, \cite{Ostrovsky3} etc. \par
 Moreover, the spaces $ G(\cdot) $ are rearrangement invariant (r.i.) spaces
with the fundamental function $ \phi(G(\psi), \delta) =
\phi(\delta) = $

$$
 \sup \{||I(A)||G, \ A \in \Sigma, \ \mu(A)
\le \delta \}, \ \delta \in (0,\infty).
$$
 If the measure $ \mu $ is nonatomic, $ \phi(G(\psi,\delta) = ||I(A)||G(\psi),$
were
$ \ \mu(A) = \delta, \ $ we have for the spaces $ G(\psi), \
\psi(\cdot) \in U \Psi,\ B( \psi) = (a,b), \ b \le \infty $
$$
\phi(G(\psi), \delta) = \sup_{p \in (a,b)} \left[ \delta^{1/p} /\psi(p) \right]. \eqno(1.5)
$$

As a slight consequence (for nonatomic measures):

$$
\phi(G(\psi), 0+) = 0; \ \lim_{\delta \to 0+} \delta/\phi(G(\psi), \delta)
= 0.
$$

 Note that in the case $ b < \infty $
$$
\delta \le 1 \ \Rightarrow C_1 \delta^{1/a} \le \phi(G,\delta) \le
C_2 \delta^{1/b},
$$

$$
\delta > 1 \ \Rightarrow C_3 \delta^{1/b} \le \phi(G,\delta) \le
C_4 \delta^{1/a}.
$$
 For instance, define in the case $ b < \infty \ \delta_1 =
\exp(\alpha h^2/(h-a)), \ \delta \ge \delta_1 \ \Rightarrow $
$$
p_1 =p_1(\delta) = \log \delta/(2 \alpha) -\left[0.25
\alpha^{-2}\log^2 \delta -a \alpha^{-1}\log \delta \right]^{1/2},
$$
$$
\phi_1(\delta) = \delta^{1/p_1}(p_1-a)^{\alpha};
$$
$$
\delta \in (0,\delta_1) \ \Rightarrow \phi_1(\delta) =
\delta^{1/h}(h-a)^{\alpha};
$$
$$
\delta_2 = \exp(-h^2\beta/(b-h)), \ \delta \in (0,\delta_2) \
\Rightarrow
$$
$$
p_2 = p_2(\delta) = - |\log \delta|/2 \beta + \left[ \log^2
(\delta/(4 \beta^2)) + b |\log \delta|/ \beta \right]^{1/2},
$$
$$
\phi_2(\delta) = \delta^{1/p_2(\delta)} (b-p_2(\delta))^{\beta};
$$
$$
\delta \ge \delta_2 \ \Rightarrow \phi_2(\delta) = \delta^{1/h}
(b-h)^{\beta}.
$$
 We obtain after some calculations:
$$
b < \infty \ \Rightarrow \phi(G(a,b;\alpha,\beta),\delta) = \max
\left[ \phi_1(\delta), \phi_2(\delta) \right].  \eqno(1.6)
$$
 Note that as $ \delta \to 0+ $

$$
\phi(G(a,b,\alpha,\beta), \delta) \sim ( \beta b^2/e)^{\beta} \
\delta^{1/b} \ |\log \delta|^{-\beta},
$$
and as $ \delta \to \infty $

$$
\phi(G(a,b,\alpha,\beta),\delta) \sim (a^2 \alpha/e)^{\alpha}
\delta^{1/a} \ (\log \delta) ^{-\alpha}.   \eqno(1.7)
$$

 In the case $ b = \infty, \beta < 0 $ we have denoting
$$
\phi_3(\delta) = ( \beta /e)^{\beta} \ |\log \delta|^{-|\beta|},
\ \delta \in (0,\exp(-h |\beta|)),
$$

$$
\phi_3(\delta) = h^{-|\beta|} \delta^{1/h}, \ \delta \ge \exp(-h
|\beta|):
$$

$$
\phi(G(a,\infty; \alpha, - \beta), \delta) = \max(\phi_1(\delta),
\phi_3 (\delta)), \eqno(1.8)
$$
and we receive as $ \delta \to 0+ $ and as $ \delta \to \infty $
correspondingly:

$$
\phi(G(a,\infty;\alpha, - \beta), \delta) \sim (\beta)^{|\beta|}
|\log \delta|^{-|\beta|},
$$
$$
\phi(G(a,\infty; \alpha,-\beta),\delta) \sim ( a^2
\alpha/e)^{\alpha} \ \delta^{1/a} (\log \delta)^{-a}.
$$
\vspace{2mm}

{\sc Boyd's indices.}

\vspace{2mm}

 At the end of this section we give using this results the expression for
the so-called Boyd's (and other) indices of $ G(\psi,a,b) $ spaces in the case $ X = [0, \infty) $ with usually Lebesgue measure. This indices play a very
important role in the theory of operators interpolation, theory of Fourier series in the r.i. spaces etc.; see  \cite{Bennet1}, \cite{Krein1}.\par

 Recall the definitions. Introduce the (linear) operators
$$
\sigma_s f(x) = f(x/s), \ s > 0,
$$
then for arbitrary r.i. space $ G $ on the set $ X = R^1_+ $

$$
\gamma_1(G) = \lim_{s \to 0+} \log ||\sigma_s||/\log s;
$$
$$
\gamma_2(G) = \lim_{s \to \infty} \log ||\sigma_s||/\log s;
$$

  We obtained (see \cite{Ostrovsky5})
$$
\gamma_1(G(\psi,a,b))= 1/b, \ \gamma_2(G(\psi,a,b))= 1/a. \eqno(1.10)
$$

 Note that also
$$
\gamma_1(G^o(\psi,a,b))= 1/b, \ \gamma_2(G^o(\psi,a,b))= 1/a.  \eqno(1.11)
$$

\vspace{2mm}

{\bf 2. \ Associate and dual spaces.} \\

\vspace{2mm}

 The complete description of the spaces conjugated (or, on the other words,
dual) to the (linear topological) spaces $ \cap_p L_p $ or conversely
$ \cup_p L_p $ see in  \cite{Davis1}, \cite{Steigenwalt1}. \par
 The spaces which are conjugate to the
 Orlicz's spaces are described, e.g., in \cite{Rao2}, chapter \
3; \cite{Rao1},pp. 123 \ - \ 142. \par
 It is easy to verify using the well - known theorem of Radon - Nicodim
that the structure of linear continuous functionals $ l = l(f) $
over the space $ G^0(\psi) = GA = GB $ is follow: $ \forall l \in
G^{o*}(\psi) (= GA^*(\psi) = GB^*(\psi)) \ \Rightarrow \exists g:
X \to R, $
$$
 l(f) = l_g(f) = \int_X f(x) g(x) \ \mu(dx), \eqno(2.0)
$$
where $ g $ is some local integrable function:

$$
\forall A \in \Sigma, \mu(A) \in (0,\infty)  \ \Rightarrow \int_A |g| \ d \mu
< \infty.
$$

 We will call as usually the space of all {\it continuous} \ in $ G(\psi, a,b)$
space functionals of a view (2.0) as {\it associated space} and
will denote as $ G^/(\psi) = G^/(\psi, a,b). $ It is evident that
$ G^/(\psi) $ is closed subspace of $ G^*(\psi). $ \par

\vspace{2mm}

{\bf Definition 3.}\par

\vspace{2mm}

 Analogously to the works \cite{Capone1}, \cite{Fiorenza1}, \cite{Fiorenza2},
\cite{Fiorenza3} we can describe the associated spaces to
the $ G(\psi, a,b) $ spaces.\par

 Let us introduce for the set $ A = (a,b) \subset[1, \infty) $ its adjoin set:

$$
 A^/ \stackrel{def}{=} \left(b^/, a^/ \right), \ a^/ = a/(a-1), b^/ = b/(b-1).
$$

\vspace{2mm}

{\it We denote by the symbol $DSL = DSL(\psi) = DSL(\Psi, a,b) = DS $
 a Banach space of   all the \ measurable functions $ g:X \to R^1 $
 with finite norm }
$$
||g||DSL(\psi) \stackrel{def}{=} \inf \left\{\sum_{k = 1}^{\infty}\psi(q(k)) \
|g_k|_{q(k)^/} \right\}, \eqno(2.1)
$$
 where $ "\inf" $ is calculated over all the sequences of measurable functions
$ g_k = g_k(x) $ such that
$$
\sum_k g_k(x) = g(x) \eqno(2.2)
$$
a.e., $ g_k \in L_{q(k)^/}, $ and all sequences of numbers $ \{q(k) \} $
belonging to the set $ (a,b): q(k) \in (a,b). $ \par

\vspace{2mm}

{\bf Remark 2.1} It may be proved analogously \cite{Capone1}, \cite{Fiorenza1},
\cite{Fiorenza3} etc. that the $ DSL(\psi) $ spaces relative the introduced norm
are really Banach function spaces and, moreover, are r.i. spaces.\par

\vspace{2mm}

{\bf Remark 2.2} It is easy to verify that in the decomposition (2.2) in the
case if $ g(x) \ge 0 $ (almost everywhere)  all the
functions $ g_k(\cdot) $ may be choose  to be non-negative.\par

{\bf Remark 2.3}.
The sequence of the numbers $ \{ q(k) \} $ in the definition (2.1)
may be choose such that

$$
\underline{\lim}_{k \to \infty} q(k) = a, \ q(k) > a,
$$
and

$$
\overline{\lim}_{k \to \infty} q(k) = b, \ q(k) < b.
$$

\vspace{2mm}

{\bf Theorem 2.1.} Let $ \psi \in U \Psi. $  Let also the triplet $ (X, \Sigma,
\mu) $ be resonant. \ Then  the spaces $ DSL(\psi) $ and $ SL(\psi) $ coincides
(up to norm equality):

$$
DSL(\psi) =  SL(\psi).  \eqno(2.3)
$$

 On the other words, the spaces $ DSL(\psi) $ are  associate to the Bilateral
Grand Lebesgue spaces $ G(\psi). $ \par

{\bf Proof.}
 This proposition  may be proved analogously to the case one-side Grand
Lebesgue spaces considered in \cite{Capone1}, \cite{Fiorenza1}-\cite{Fiorenza3}. \par

\vspace{2mm}

{\bf Step 1}. We prove at first the implication $ g \in
SL(\psi) \ \Rightarrow l_g(\cdot) \in GL^/(\psi), $ with norm equality. \par

  Namely, if $ g = \sum_k g_k, \ f \in G(\psi), \  ||f||G(\psi) = 1 $ and for some sequence $ \{q(k)\} \in (a,b) $ and for any $ \epsilon \in (0,1) $

$$
\sum_k \psi(q(k)) |g_k|_{(q(k))^/} \le ||g||GS(\psi) + \epsilon= C
< \infty,
$$
then $ |f|_{q(k)} \le \psi(q(k)), \ q(k) \in (a,b). $ We have
using H\"older inequality:

$$
|l_g(f)| \le \sum_k |f|_{q(k)} \cdot |g_k|_{(q(k))^/} \le \sum_k
\psi((q(k))) \cdot |g_k|_{(q(k))^/} = C < \infty,
$$
or, equally,

$$
 |l_g(f)| = \left|\int_X f(x) \ g(x) \ \mu(dx) \right|
\le ||f||G(\psi) \ ||g||SL(\psi). \eqno(2.4)
$$
 Note that the inequality (2.4) is called as usually as generalized H\"older
inequality. \par

\vspace{2mm}

{\bf Step 2}. Now we prove the inverse implication. \par
It is sufficient to consider only the case of the space $ G^o(\psi; a,b). $\par

 Let $ f(\cdot) $ be arbitrary element of the space $ G^o(\psi) $ with unit norm:

$$
|f|_p \le \psi(p), \ p \in (a,b);
$$
and let $ l_g $
be  a linear functional on the space $ G^o(\psi) $ of the standard view:

$$
l_g(f) = \int_X f(x) \ g(x) \ d \mu.
$$
 We deduce using the classical H\"older inequality:

$$
\forall p \in (a,b) \ |l_g(f)| \le |f|_p \ |g|_{p^/} \le \psi(p) \  |g|_{p^/}.
$$

Therefore,

$$
||g||SL(\psi) \le \inf_{q \in A^{/} } \psi(q^{/} ) |g|_q. \eqno(2.5)
$$

 It remain to prove  the finiteness of the right side of inequality (2.5), or
equally the following implication:

$$
\exists q \in A^/, \ |g|_q < \infty. \eqno(2.6)
$$
We will prove  this proposition by the method "reduction ad absurdum". Let us
consider for definiteness the case of purely atomic measure $ \mu:  X =
(1,2,\ldots;)  \ \mu({k}) = 1. $ \par
Then the linear functional $ l_g $ has a view:

$$
l_g(f)  = \sum_{n=1}^{\infty} f_n \ g_n
$$
and we suppose
$$
\forall q \in A^/ \ \sum_{n=1}^{\infty} (|g|_n)^q = \infty. \eqno(2.7)
$$
We will suppose without loss of generality  that $ f_n \ge 0, \ g_n \ge 0, \
n = 1,2, \ldots. $ \par

 Define the truncated sequence

$$
f^{(N)}(n) = f_n, \ n \le N, \ f^{(N)}(n) = 0, \ n > N;
$$
and the correspondent truncated sum
$$
S_N(\epsilon) = \sum_{n=1}^N |q_n|^{a^/ - \epsilon}, \eqno(2.8)
$$
$ \epsilon = \const \in (0, \epsilon_0).$ \par
 It follows from the uniform boundedness principle that there exists the positive finite constant $ K, \ K \in (0,\infty), $ such that

$$
L:= \sum_{n=1}^{\infty} f_n \ g_n \le K \cdot \sup_{p \in (a,b)}
\frac{|f|_p}{\psi(p)}. \eqno(2.9)
$$
 We choose in the inequality (2.9) the sequence $ f $ as follows:

$$
f_n = (g_n)^{a^/-(\epsilon + 1) }, \ n \le N
$$
and $ f_n = 0 $ in other case.\par

Substituting into (2.9), we get:

$$
L \le K \cdot \sup_{p \in (a,b)} \frac{ (S_N(\epsilon))^{1/p}}{\psi(p)}=
K \phi(G(\psi), S_N(\epsilon)), \eqno(2.11)
$$
where, recall, $ \phi(\cdot, \cdot) $ denotes the fundamental function. \par

 We obtain from (2.11), denoting $ X = S_N(\epsilon): $

$$
X \le K \cdot \phi(G(\psi),X). \eqno(2.12)	
$$

 The solving of inequality (2.12) has a view: $ X \le C, $ where the constant
$ C $ does not depend on the $ \epsilon $ and $ N. $ Therefore,

$$
\sup_{\epsilon} \sup_N S_N(\epsilon) \le C,
$$
or equally

$$
\sum_{n=1}^{\infty} |g_n|^{a^/} < \infty
$$
and  hence  for all the values $ q \in (a^/, b^/) $

$$
\sum_{n=1}^{\infty} |g_n|^{q} < \infty,
$$
in contradiction.\par
This completes the proof of theorem 2.1.\par

\vspace{2mm}

{\bf Step 3}.\par

 Note analogously to the article \cite{Fiorenza1} that the inequality
(2.5) is exact.\par
We intend to prove that for all $ f \in G^o(\psi), \ f \ge 0 $ there exists
a function $ g \in \cup_{ q \in (b^/,a^/)} L_q $ for which
$$
|l_g(f)| = ||f||G(\psi) \ \cdot \ ||g||SL(\psi).
$$
 Indeed, the function $ h(p) = |f|_p/\psi(p), p \in [a,b] $
is continuous and $ h(a) = h(b) = 0; $ therefore there exists a
value $ \sigma \in (a,b) $ for which $ |f|_{\sigma} = ||f|| \cdot
\psi(\sigma). $ We have choosing the function $ g $ such that

$$
l_g(f)= |f|_{\sigma} \ |g|_{\sigma^/}:
$$

$$
||f||G(\psi) \ ||g||GS(\psi) \ge l_g(f) =
$$

$$
\psi(\sigma) ||f|| \ |g|_{\sigma^/} = ||f|| \ \left[|g|_{\sigma^/}
\psi((\sigma^/)^/) \right] \ge ||f||G(\psi) \ ||g||GS(\psi).
$$

\vspace{2mm}

{\bf Remark 2.4.}
 As a corollary: if (in addition) our measure $ \mu $ is separable, then
we infer:

$$
G(\psi)^/ = GB(\psi)^* = GA(\psi)^* = G_0(\psi)^*,
$$

see \cite{Bennet1}, p. 20 \ - \ 22. \par

\vspace{3mm}

{\bf 3. \ Norm's absolutely continuity. Compactness.} \\

\vspace{2mm}

 We will say as usually, see \cite{Bennet1}, pp. 14-16
that the function $ f $ belongs to some r.i. space $ Y $  over source triplet
$ (X,\Sigma, \mu) $ with the norm $ ||\cdot||Y $ has
{\it absolutely continuous norm in this space and write }
$ f \in YA, $  {\it if}

$$
 \lim_{\delta \to 0} \sup_{A: \mu(A) \le \delta} ||f \ I_A||Y=0. \eqno(3.1)
$$

We will write in the case $ Y = G(\psi): \ G(\psi)A = GA(\psi), \ G(\psi)B = GB(\psi) $  and analogously $ SL(\psi)B = SLB(\psi), SL(\psi)A =
SLA(\psi). $\par

\vspace{3mm}

{\bf Theorem 3.1.} {\it The space $ SL(\psi; a,b) $ satisfies the
absolutely continuity norm property}.\par

\vspace{3mm}

{\bf Proof.} The proof is at the same as in \cite{Fiorenza1}. Namely,
let $ g \in SL(\psi), \ ||g||SL(\psi) = 1. $  We write the  following
consequence from the decomposition for this function:

$$
\sum_{k=1}^{\infty}  \left[ \psi(q(k)) \ |g_k|_{q(k)^{/} }  \right] \le 2.
$$
 Let also $ \{ E(n) \}, n = 1,2, \ldots $ be monotonically decreasing sequence
of measurable sets  such that  as $ n \to \infty $

$$
\mu(E(n)) \downarrow 0.
$$
 We can suppose without loss of generality  $ \mu(E(1)) \le 1.$
 Let us denote as in \cite{Fiorenza1}

 $$
 a(k,n) = \psi(q(k)) |g_{k} | \ I(E(n))|_{q(k)^{/} }.
 $$
We have  for all values $ n: $

$$
||g \ I(E(n))||SL(\psi) \le \sum_{k=1}^{\infty} a(k,n) < \infty.
$$

 Since

 $$
 \sum_{k=1}^{\infty} a(k,n) < \infty
 $$
and $ \forall k \ \Rightarrow a(k,n) \downarrow 0 $ as $  n \uparrow \infty, $
we deduce:

$$
\lim_{n \to \infty} ||g \ I(E(n))||SL(\psi) = 0, \eqno(3.2)
$$
QED.\par

\vspace{3mm}

{\bf Consequences.} \par

\vspace{2mm}
 The next assertions follows from the theorem 3.1 and from the well \ - \ known
facts about the general theory of Banach functional  and rearrangement invariant
spaces (\cite{Bennet1}, chapters 1 \ - \ 2). \par

 In the Bilateral Small Lebesgue spaces are true the Levi's
theorem of a monotone convergence, Fatou property
 and Lebesgue majoring convergence theorem.\par
 Both the spaces $ G(\psi) $ and $ SL(\psi) $ are the interpolation spaces
between the spaces $ L_1(X,\mu) $ and $ L_{\infty}(X,\mu) $ relative the
real method of interpolation. \par
 If we assume in addition to  the conditions of
theorem 3.1 that the measure $ \mu $ is nonatomic and separable,
 then the space $ BSL(\psi) $ is also separable, coincides with $ BSLA(\psi), BSLB(\psi) $ and moreover

 $$
[BSL(\psi)]^* = [BSL(\psi)]^/ = BGL(\psi). \eqno(3.3)
 $$
but they are not reflexive. \par

 Note for comparison:
the subspaces $ BGA(\psi), BGB(\psi), BG^0(\psi) $ are closed {\it strictly}
subspaces of the space $ BG(\psi). $  The following important property of the Grand Lebesgue spaces is proved in \cite{Ostrovsky1}, \cite{Ostrovsky2}.

 {\it If } $ \psi \in U\Psi, $ {\it then }

$$
G^0(\psi) = GB(\psi) = GA(\psi).
$$

\vspace{3mm}

{\bf 4. Fundamental functions and Boyd's indices.}\par

\vspace{2mm}

We suppose in this section that the measure $ \mu $ again is sigma-finite and nonatomic.\par
The next fact follows from the well-known result about the connection
between the fundamental functions of  r.i. spaces and its associate (see,
for example, \cite{Bennet1}, chapter 2, theorem 5.2.) \par

\vspace{2mm}

{\bf Theorem 4.1.}

\vspace{2mm}

$$
\phi(SL(\psi), \delta) = \frac{\delta}{\phi(G(\psi), \delta)}, \ \delta \in
(0, \mu(X)). \eqno(4.1)
$$

We will use further for brevity
 the notation $ \chi(\delta) = \phi( SL(\psi, \delta) ). $ \par
 The assertion of the theorem 4.1 allows  us to calculate the {\it exact
 value} of  $ SL(\psi) $ norm of the indicator function. \par

Notice that

$$
\chi(0+) = 0; \ \lim_{\delta \to \infty} \chi(\delta) = \infty.
$$

\vspace{2mm}

{\bf Lemma 4.1.} Let $ A $ be a measurable set: $ A \in \Sigma $ such that
$ 0 < \mu(A) < \infty. $  The following equality holds:

\vspace{2mm}

$$
|| \ I(A) \ ||SL(\psi) = \frac{\mu(A)}{\phi(G(\psi), \mu(A))} = \chi(\mu(A)).
\eqno(4.2)
$$

 Recall that the exact values for the fundamental functions for some $ G(\psi) $
spaces there are in the formulas (1.6), (1.8).  \par
 Let us give an other example for $ SL(\psi) $ norm exact value calculation. \par
{\bf Example 4.1.}  Assume that on the source triplet $ (X,\Sigma, \mu) $
and  assume that for some interval $ (a,b), \ 1 \le a < b \le \infty $  there exists a  value $ q_0, a < q_0 < b $ and a measurable non \ - \ negative
function $ g: X \to R $ for which

$$
|g|_q < \infty \Leftrightarrow q = q_0.
$$

 This condition is satisfied, for instance, when $ X = R^d $ and $ \mu $ is
$ d \ - \ $ dimensional Lebesgue measure. \par
We deduce using the explicit expression for the $ SL(\psi) $ norm  that for
this function $ g $ the inequality (2.5) transforms to the equality:

$$
||g||SL(\psi) =  \psi \left(q_0^{/} \right) |g|_{q_0}. \eqno(4.3)
$$

 The assertions of the theorem 4.1 and lemma 4.1 give us the possibility to
compute the important for applications  Boyd's $ \gamma_1, \ \gamma_2 $
and Shimogaki's $ \beta_1, \ \beta_2 $ indices
for Bilateral Small Lebesgue spaces on the basis of ones results for Grand Spaces,
which are calculated in \cite{Ostrovsky3}, \cite{Ostrovsky5}. Indeed: \par

\vspace{2mm}

{\bf Lemma 4.2.} Let $ \psi \in \Psi(a,b), \ 1 \le a < b \le \infty. $ Then:

$$
{\bf A.} \gamma_1(SL(\psi)) = 1-1/a; \ \gamma_2(SL(\psi)) = 1-1/b;  \eqno(4.4)
$$

$$
{\bf B.} \beta_1(SL(\psi)) = 1-1/a; \ \beta_2(SL(\psi)) = 1-1/b;  \eqno(4.5)
$$

{\bf C.} The Hardy \ - \ Littlewood maximal operator $ M $
in the case $ X = R^d $
is in bounded in the Small Lebesgue space $ SL(\psi). $ \par
{\bf D.} The Hilbert's transform $ H $ in the case $ X = R^d $
is in bounded in the Small Lebesgue space $ SL(\psi) $ iff $ a > 1. $ \par
{\bf E.} Let here $ X = [0, 2 \pi). $  The Fourier series for arbitrary function
$ f \in SL(\psi) $ convergence to $ f $ in the $ SL(\psi) $ norm iff
 $ a > 0, \ b < \infty. $\par

\vspace{3mm}

{\bf 5. \ Improving of the norm estimation.} \\

\vspace{2mm}

Recall that

$$
||I(A)||SL(\psi) = \chi(\mu(A)). \eqno(5.1)
$$

The equality (5.1) gives us the possibility to {\it estimate} the $ SL(\psi) $ norm
of arbitrary function belonging to the space $ SL(\psi). $ \par
 Consider a {\it simple} function of a view:

 $$
 f(x) = \sum_{k=1}^n c(k) \ I(H_k,x), n < \infty, \eqno(5.2)
 $$
$ H_k \in \Sigma, \ k \ne l \ \Rightarrow H_k \cap H_l = \empty, \ c(k) \in R; $
i.e. $ f \in SLB(\psi) $  and conversely the set of all the functions of a view (5.2)
is dense in the $ SLB(\psi) $ space and following in all space $ SL(\psi). $ \par
 We will call the decomposition (5.2) as a representation of a function $ f $
 and will denote the set of all simple functions as $ Sim(\psi). $ \par

 We have for these functions:

 $$
 ||f||SL(\psi) \le \sum_{k=1}^n |c(k)| \ \chi(\mu(H_k)) \stackrel{def}{=}
 \int_X  f(x) \ d \chi. \eqno(5.3)
 $$

 Let now $ f(\cdot) $ be  arbitrary non \ - \ negative function from the space
 $ SL(\psi). $  Denote by $ R(f) $ the set of  all the simple functions $ \{ g \} $
 greatest than $ f: g \ge f. $  We define:

 $$
 \int_X f(x) \  d \chi \stackrel{def}{=} \inf_{g \in R(f)}  \int_X g(x) \ d \chi.
 \eqno(5.4)
 $$

 For the function $ f \in SL(\psi) $ which may not to be  non-negative we define

 $$
 \int_X f(x) \ d \chi =  \int_X f_1(x) \ - \ \int_X f_2(x) \ d \chi,
 \eqno(5.5),
 $$
where the non \ - \ negative functions $ f_1 $ and $ f_2 $ in (5.5) from the set
$ Sim(SL) $ are such that the function $ f $ is the difference  of the functions
$ f_1 $ and $ f_2: \ f(x) = f_1(x)-f_2(x) $ and

$$
\int_X [ |f_1(x)| + |f_2(x)|] \ d \chi = \inf \left\{ \int_X [ |g_1(x)| + |g_2(x)|] \
d \chi \right\},
$$
where
 $ g_1, g_2 \in Sim(SL), \ g_1, g_2 \ge 0 $ and $ g_1(x)-g_2(x) = f(x). $ \par

For the unbounded functions $ f = f(x) $ we may define

$$
\int_X f(x) \ d \ \chi = \lim_{N \to \infty} \int_X [ f(x) \ I(|f(x)| \le N ] \ d \chi,
$$
if there exists.\par
 Let us denote
$$
|||f||| = |||f|||SL(\psi) = \int_X |f(x)| \ d \chi. \eqno(5.6)
$$

 Note that the functional $ f \to |||f|||SL(\psi), \ f \in SL(\psi) $
 obeys the following norm properties:

 $$
|||f||| \ge 0; \ |||f||| = 0 \Leftrightarrow  f = 0 \ (\mod \ \mu);
 $$

$$
||| \lambda \ f||| = |\lambda| \cdot |||f|||, \ \lambda \in R;
$$

$$
||| f + g ||| \le |||f||| + |||g|||.
$$

 Since the set of the simple functions is dense in the space $ SL(\psi), $ we
 obtain the following convenient for applications estimation.\par

 {\bf Theorem 5.1.}

 $$
 ||f||SL(\psi) \le |||f|||SL(\psi). \eqno(5.7)
 $$

{\bf Remark 5.1.} Note that still in the case $ \mu(X) = 1 $ the functional $ f \to
|||f|||SL(\chi) $ is discontinued in the {\it uniform} norm
$ |f|_{\infty} = \vraisup_{x \in X} |f(x)| $ in all the points aside from the origin.
Let us consider the following  example.\par
{\bf Example 5.1.} Let $ X = [0,1] $ and $ \ n = 2,3,4,\ldots. $ Define the functions
$ f(x) = 1 $ and the sequence of a functions

$$
f_n(x) = I( x  \in [0,1/2] ) + (1 + 1/n) I(x \in (1/2, 1)). \eqno(5.8)
$$
Let also $ \chi(\delta) = \sqrt{\delta}, \ \delta \in (0,1). $ We observe:

$$
\lim_{n \to \infty} |f_n-f|_{\infty} = 0; \ |f|_{\infty} = 1;
$$
but
$$
|||f_n|||SL(\chi) =  \sqrt{0.5} + (1 + 1/n) \sqrt{0.5} \to \sqrt{2}, \ n \to \infty.
$$

\vspace{3mm}

{\bf 6. \ Convergence and compactness. } \\

\vspace{2mm}

 {\bf A.} Note that in the case if $ X $ is the convex bounded subset of $ R^n $ with
usually {\it bounded} Lebesgue measure, for the spaces $ G^o(\psi) $  and
$ SL(\psi) $   are true
the classical conditions of Riesz's and Kolmogorov's for compactness
of some subset $ F = \{f_{\alpha} \} \subset G^o(\psi), $ as long as
these spaces are separable. \par

\vspace{2mm}
 {\bf B.} If some subset $ F = \{f_{\alpha} \} \subset G^o(\psi) $ is
closed, bounded, is compact set in the sense of $ L_p, p \in (a,b) $
convergence and has uniform absolutely continuous norm: $ F \in
UCN, $ then $ F $ is compact set in the space $ G(\psi). $ \par

 \vspace{2mm}
 {\bf C.} The direct calculation of the norm $ ||f||SL(\psi) $ is very hard. But
 we can replace the distance $ || f-g ||SL(\psi) $ by the distance

 $$
 d_{\psi}(f,g) \stackrel{def}{=} |||f-g|||SL(\psi) \eqno(6.1)
 $$
 in order to formulate the {\it sufficient } conditions for convergence and
 compactness in the Bilateral Small Lebesgue spaces. For instance, if $ f_n, f
 \in SL(\psi) $ and $ |||f_n-f||| \to 0, \ n \to \infty, $ then
  $ ||f_n-f||SL(\psi) \to 0, \ n \to \infty. $\par

  If some closed subset $ U $ of the space $ SL(\psi) $ is compact set in the
  distance $ d_{\psi}, $ then $ U $ is compact set relative the source norm
  of the space $ SL(\psi) $ etc. \par

 \vspace{4mm}



\begin{thebibliography}{69}

\bibitem{Bennet1}	

 C. Bennet, R. Sharpley. {\it Interpolation of operators.} Orlando,
Academic press Inc., (1988). \\

\bibitem{Capone1}

 C. Capone  and A. Fiorenza. {\it On Small Lebesgue spaces.} Instituto per le Applicazioni del Calcolo "Mauro  Picone", Napoli,  (2004), Rapporto tecnico n. 278/04.\\

\bibitem{Carro1}

M. Carro, J. Martin. {\it Extrapolation theory for the real interpolation method.} Collect. Math. 33 (2002), 163-186.\\

\bibitem{Davis1}

 H.W. Davis, F.J.Murray, J.K.Weber. {\it Families of $ L_p - $
spaces with inductive and projective topologies.} Pacific J.Math.
v. 34, (1970), 619-638.\\

\bibitem{Fiorenza1}

 A.Fiorenza. {\it Duality and reflexivity in grand Lebesgue
spaces.}
Collectanea Mathematica (electronic version), {\bf 51}, 2, (2000), 131-148.\\

\bibitem{Fiorenza2}

 A. Fiorenza and G.E. Karadzhov. {\it Grand and small Lebesgue
spaces and their analogs.} Consiglio Nationale Delle Ricerche,
Instituto per le Applicazioni del Calcoto Mauro Picine", Sezione
di Napoli, Rapporto tecnico n. 272/03, (2005).\\

\bibitem{Fiorenza3}

A.Fiorenza, J.-M. Rakotoson. {\it On Small Lebesgue spaces and their
applications.}
Comptes  Rendus Mathematique, Volume  334, Number 1, 1 January
(2002), pp. 23-26. \

\bibitem{Fiorenza4}

 A.Fiorenza, Rakotoson J.-M. {\it New properties of small Lebesgue
spaces and their applications.}
Comptes  Rendus Mathematique, Volume  349, Number 3, 22 April
(2008), pp. 147-181. \\

\bibitem{Iwaniec1}

 T.Iwaniec and C. Sbordone.{\it On the integrability of the Jacobian under minimal hypotheses.} Arch. Rat.Mech. Anal., 119, (1992), 129-143.\\

\bibitem{Iwaniec2}

T.Iwaniec, P. Koskela and J. Onninen. {\it Mapping of finite
distortion: Monotonicity and Continuity,} Invent. Math. 144 (2001), 507-531.\\

\bibitem{Jawerth1}

B. Jawerth, M.Milman. {\it Extrapolation theory with applications.} Mem. Amer. Math. Soc. 440 (1991).\\

\bibitem{Karadzhov1}

G.E. Karadzhov, M. Milman M. {\it Extrapolation theory: new
results and applications. } J. Approx. Theory, 113 (2005), 38-99.\\

\bibitem{Krasnoselsky1}

M.A.Krasnoselsky, Ya.B.Rutisky. {\it Convex functions and
Orlicz's Spaces.} P. Noordhoff LTD, The Netherland, Groningen, 1961. \\

\bibitem{Kozatchenko1}

Yu.V. Kozatchenko,  E.I. Ostrovsky. {\it Banach spaces of random
variables of subgaussian type.} Theory Probab. And Math. Stat.,
Kiev, (1985), p. 42-56, (in Russian).\\

\bibitem{Krein1}

S.G.Krein, Yu. Petunin, and E.M.Semenov. {\it Interpolation of
linear operators.} New York, AMS, (1982);\\

\bibitem{Ostrovsky1}

 E.Ostrovsky. {\it Exponential Orlicz's spaces: new norms and
applications.}
Electronic Publications, arXiv/FA/0406534, v.1, (25.06.2004.)\\

\bibitem{Ostrovsky2}

 E. Ostrovsky.{\it Exponential Estimations for Random Fields.}
Moscow - Obninsk, OINPE, (1999) (in Russian).\\

\bibitem{Ostrovsky3}

 E.Ostrovsky, L.Sirota. {\it Some new rearrangement invariant
spaces: theory and applications.} Electronoc publications:
arXiv:math.FA/0605732 v1,29, (May 2006);\\

\bibitem{Ostrovsky4}

 E.Ostrovsky, L.Sirota. {\it Fourier Transforms in Exponential
Rearrangement Invariant Spaces.} Electronoc publications:
arXiv:math.FA/040639, v1, (20.6.2004.)\\


\bibitem{Ostrovsky5}

E. Ostrovsky, L. Sirota.
{\it Moment Banach spaces: Theory and applications.} HAIT Journal of Science
and Engineering, C, Holon, ISRAEL,(2007), V. 4, Issues 1-2, pp. 233-262.


\bibitem{Rao1}

 M.M. Rao, Z.D.Ren. {\it Theory of Orlicz Spaces.} Basel-New
York, Marcel Decker, (1991);\\

\bibitem{Rao2}

 M.M. Rao, Z.D.Ren. {\it Application of Orlicz Spaces.} Basel-New York, Marcel Decker, (2002. \\

\bibitem{Rao3}

 M.M.Rao. {\it Measure Theory and Integration.} Basel-New
York, John Wiley, Marcel Decker, second Edition, (2004);\\

\bibitem{Seneta1}

 E. Seneta. {\it Regularly Varying Functions}. Mir, Moscow
edition, (1985).\\

\bibitem{Steigenwalt1}
 M.S.Steigenwalt and A.J.While. {\it Some function spaces
related to} $ L_p. $ Proc. London Math. Soc, 22, (1971), 137-163;\\


\end{thebibliography}
\end{document}